\documentclass{tsp-short}
\usepackage{graphicx} 
\usepackage{geometry}
\usepackage{hyphenat}
\usepackage[T2A]{fontenc}
\usepackage[english]{babel}
\usepackage{tikz}
\usetikzlibrary{positioning,chains,fit,shapes,calc,arrows,patterns,external,shapes.callouts,graphs,hobby,intersections,spath3}
\usepackage{amsfonts}
\usepackage{amsmath}
\usepackage{amssymb}
\usepackage{bbm}
\usepackage[utf8]{inputenc}
\usepackage{amsthm}
\usepackage{array}
\usepackage{mathrsfs}
\usepackage[skip=10pt plus1pt, indent=0pt]{parskip}
\usepackage{cancel}
\usetikzlibrary{positioning}
\usetikzlibrary{petri}
\tikzset{state/.style={inner sep=0pt,minimum size=7mm,label=center:$#1$,name=#1},
	redarrow/.style={->, black, fill=none},bluearrow/.style={->, black, fill=none},  
	redline/.style={-,black,fill=none, thick},blueline/.style={<->,black,fill=none}}
\usepackage[backend=biber]{biblatex}
\usepackage{csquotes}
\usepackage{braids}
\usepackage[all,knot]{xy}
\usepackage{tabularray}
\newcolumntype{C}[1]{>{\centering\let\newline\\\arraybackslash\hspace{0pt}}m{#1}}
\usepackage{mathtools}
\usepackage{amssymb}
\usepackage{tikz}
\usetikzlibrary{shapes,snakes}
\usepackage{adjustbox}
\usepackage{caption}

\DeclareMathOperator{\EX}{\mathbb{E}}

\newtheorem{thm}{Theorem}[section]

\newtheorem{cor}{Corollary}[section]
\newtheorem{prp}{Proposition}[section]
\theoremstyle{remark}

\theoremstyle{definition}

\begin{document}
	
	\title
	[Random braids and random walks on finite groups]
	{Random braids and random walks on finite groups}
	
	\author{H.M. Zhylinskyi}
	\address{Department of Mathematics, Jagiellonian University, Krakow,
		Poland}
	\email{georgyzhilinsky144@gmail.com}
	
	
	\begin{abstract}
		In this paper, we study properties of random walks on finite groups and later use them to obtain the limiting braid length expectation and component number of braid closure in a model of random braids, which is constructed by lifting elements of random walk on a Coxeter group to a braid group.
	\end{abstract}
	
	\maketitle
	
	\section{Introduction}
	
	We will consider random walks where at each step with some probability an element of a group can be multiplied by one of the generators from the presentation of this group. We will be interested in finding the limiting expectations of different functions on this group, especially the length function. In \cite{dorogovtsev2020loop}, a similar construction of loop-erased random walks on finite groups is studied where loops are erased in order they appear. Also, in \cite{eriksen2005expected} and \cite{bousquet2009expected} results of the same type as in Proposition $3.1$ are given, but for the symmetric group $S_n$, and \cite{sjostrand2012expected}, where the same problem is considered for certain Coxeter groups. Furthermore, we mention \cite{nechaev1996random}, \cite{voituriez2002random} where random walks on braid groups are also considered.
	
	\section{Braids and braid groups}
	
	Braid on $n$ strings is an object consisting of two horizontal lines $L_a$ and $L_b$ in $\mathbb{R}^3$ containing two ordered sets of points $a_1 = (1,0,0),\ldots, a_n = (n,0,0) \in L_a$ and $b_1 = (1,0,1),\ldots, b_n = (n,0,1) \in L_b$, and $n$ strands that satisfy the following properties:
	\begin{itemize}
		\item Each strand connects $a_i$ with $b_j$ for some $i,j$ and strands are pairwise disjoint.   
		\item Strands have to strictly ascend all the way up.
	\end{itemize}
	
	Instead of thinking about braids in space, we can project them on a plane passing through $L_a$ and $L_b$. Additionally, through ambient isotopy, it can be achieved that there are only finitely many double crossings, which results in a braid diagram.
	
	\begin{figure}
		\begin{minipage}[!h]{0.49\linewidth}
			\centering
			\begin{tikzpicture}[scale=.8]
				\braid[line width=1pt, nudge factor=-0.05, control factor=0.75, number of strands=5] (braid) at (0,0) s_1 s_4 s_3 s_2^{-1} s_3^{-1};
			\end{tikzpicture} 
			\caption{Geometric braid on 5 strands}
			\label{fig:edit-label}
		\end{minipage}
		\begin{minipage}[!h]{0.49\linewidth}
			\centering
			\begin{tikzpicture}[scale=.8, rotate=90]
				\braid[line width=1pt, nudge factor=-0.05, control factor=0.75] (braid) at (0,0) s_1 s_2 s_1^{-1} s_2;
				\draw[line width=1pt, dotted] (braid-rev-1-s) -- ++(0,1) -- ++(1,0) |- ([shift=(south:1)]braid-rev-3-e) -- (braid-rev-3-e);
				\draw[line width=1pt,dotted] (braid-rev-3-s) -- ++(1,1.5) -- ++(1.5,0) |- ([shift=(south:1.5)]braid-rev-2-e) -- (braid-rev-2-e);
				\draw[line width=1pt,dotted] (braid-rev-2-s) -- ++(-.5,1) -- ++(-1,0) |- ([shift=(south:1)]braid-rev-1-e) -- (braid-rev-1-e);
			\end{tikzpicture}
			\caption{Closure of a braid}
			\label{fig:edit-label}
		\end{minipage}
	\end{figure}
	
	Braid $\beta$ can be turned into a link by connecting the opposite nodes of a braid. This operation is called the closure of a braid $\beta$, and we will denote it by $\overline{\beta}$. 
	
	The braid group on $n$ strands is a group of braid equivalence classes under ambient isotopy where the composition is braid concatenation. Let $\sigma_i$ and $\sigma_i^{-1}$ be two types of crossings in the braid diagram as shown in the figure below.
	
	\begin{minipage}[!h]{0.5\textwidth}
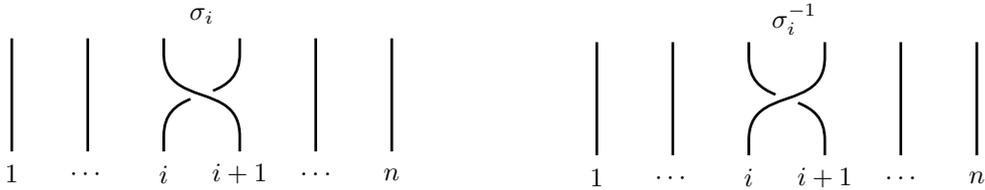
		\begin{center}
			\begin{adjustbox}{margin=0.1cm}
				\begin{tikzpicture}
					\braid[line width=1pt, nudge factor=-0.05, control factor=0.75, number of strands=6] (braid) at (0,0) s_3;
					
					\node at (0,-1.8) {$1$};
					\node at (1,-1.8) {$\cdots$};
					\node at (2.0,-1.8) {$i$};
					\node at (3,-1.8) {$i+1$};
					\node at (4,-1.8) {$\ldots$};
					\node at (5,-1.8) {$n$};
					\node at (2.5,0.3) {$\sigma_i$};
				\end{tikzpicture}
			\end{adjustbox}
		\end{center}
	\end{minipage} 
	\begin{minipage}[!h]{0.5\textwidth}
		\begin{center}
			\begin{adjustbox}{margin=0.1cm}
				\begin{tikzpicture}
					\braid[line width=1pt, nudge factor=-0.05, control factor=0.75, number of strands=6] (braid) at (0,0) s_3^{-1};
					
					\node at (0,-1.8) {$1$};
					\node at (1,-1.8) {$\cdots$};
					\node at (2.0,-1.8) {$i$};
					\node at (3,-1.8) {$i+1$};
					\node at (4,-1.8) {$\ldots$};
					\node at (5,-1.8) {$n$};
					\node at (2.6,0.3) {$\sigma_i^{-1}$};
				\end{tikzpicture}
			\end{adjustbox}
		\end{center}
	\end{minipage} 
	
	This setting was first considered by Artin and in \cite{artin1947theory} he proved the following theorem.
	
	\begin{thm}
		The braid group is isomorphic to the group with generators $\sigma_1,\ldots,\sigma_{n-1}$ subject to braid relations:
		\begin{itemize}
			\item $\sigma_i \sigma_j = \sigma_j \sigma_i$ for $|i-j| > 1$,
			\item $\sigma_i \sigma_{i+1} \sigma_i = \sigma_{i+1} \sigma_i \sigma_{i+1}$ for $i = \overline{1,n-1}$.
		\end{itemize}
	\end{thm}
	Further we will denote it by $Br_n$.
	
	\section{Random walks on finite groups}
	
	Let $S$ be a set and $\langle S \rangle$ be a free group on $S$. Define a \emph{word} on $S$ as some product of elements in $S \cup S^{-1}$. And let $R$ be a set of words on $S$. In order to form a group $G$ with presentation $\langle S \; | \; R \rangle$, take the quotient of $\langle S\rangle$ by the smallest normal subgroup $N$ containing $R$. The elements of $S$ are called \emph{generators} and the elements of $R$ are called \emph{relations} or \emph{relators}. 
	
	Define the length function $\ell_{G}(g)$ of $g \in G = \langle S \; | \; R \rangle$ to be the smallest positive integer $r$ such that there exists a word of length $r$ representing $g$. We call such an expression $g$ \emph{reduced}. By convention, $\ell(1) = 0$ where $1$ is the empty word.
	
	There is a natural way to represent a group $G$ with its generators using a \emph{Cayley graph}. Cayley graph $\Gamma = \textnormal{Cay}(G, S)$ is a directed graph $\Gamma$ that satisfies the following conditions:
	\begin{itemize}
		\item Each element of $G$ is assigned a vertex in $\Gamma$.
		\item Vertex $a$ is connected to vertex $b$ ($a \to b$) if and only if $b = a \sigma$ for some $\sigma \in S$.
	\end{itemize}
	
	It is obvious from the above definition that graph $\Gamma$ is connected, since there is a path from $1$ to each vertex. Also note that the length of an element $g \in G$ corresponds to the length of the shortest path between vertices $1$ and $g$ in $\Gamma$. 
	
	\begin{figure}[h]
		\begin{center}
			\begin{tabular}{c c}
				\begin{tikzpicture} [node distance=3cm,] 
					\graph [clockwise=6,radius=2cm] {
						v1/$(123)$,
						v2/$(213)$,
						v3/$(231)$,
						v4/$(321)$,
						v5/$(312)$,
						v6/$(132)$
					};
					\draw[<->] (v1)--(v2) node[midway,sloped,above]{$\sigma_1$};
					\draw[<->] (v2)--(v3) node[midway,sloped,above]{$\sigma_2$};
					\draw[<->] (v3)--(v4) node[midway,sloped,above]{$\sigma_1$};
					\draw[<->] (v4)--(v5) node[midway,sloped,above]{$\sigma_2$};
					\draw[<->] (v5)--(v6) node[midway,sloped,above]{$\sigma_1$};
					\draw[<->] (v6)--(v1) node[midway,sloped,above]{$\sigma_2$};
				\end{tikzpicture}
				&
				\scalebox{0.85}{
					\begin{tikzpicture}
						\node[state=1]{};
						\node[state=a,above=of 1]{};
						\node[state=a^2,right=of a]{};
						\node[state=a^3,below=of a^2]{};
						\node[state=b,below left=of 1]{};
						\node[state=ab,above left=of a]{};
						\node[state=a^2b,above right=of a^2]{};
						\node[state=a^3b,below right=of a^3]{};
						\draw[redarrow](1)--(a);\draw[redarrow](ab)--(b);
						\draw[redarrow](a)--(a^2);\draw[redarrow](a^2b)--(ab);
						\draw[redarrow](a^2)--(a^3);\draw[redarrow](a^3b)--(a^2b);
						\draw[redarrow](a^3)--(1);\draw[redarrow](b)--(a^3b);
						\draw[blueline](1)--(b);
						\draw[blueline](a)--(ab);
						\draw[blueline](a^2)--(a^2b);
						\draw[blueline](a^3)--(a^3b);
					\end{tikzpicture}
				}
			\end{tabular}
		\end{center}
		\caption{Cayley graphs of $S_3$ and $D_4$}
	\end{figure}
	
	To illustrate the above definitions, we provide an example of the Cayley graph of the symmetric group $S_3$ in Figure 3 given by the presentation $\langle \sigma_1 = (12),\;\sigma_2 = (23) \mid \sigma_1 ^2 = \sigma_2 ^ 2 = (\sigma_1 \sigma_2 )^3 = 1\rangle $ and the Cayley graph of the dihedral group $D_4$ given by the presentation $\langle a,b \mid a^4 = b^2 = 1, aba = b \rangle$.
	
	Finally, we define a random walk \textbf{RW}$(G, P)$ to be a Markov chain $X_0, X_1, \ldots$ ($Pr(X_0 = 0) = 1$) with state space $G$ and transition probability matrix $P = (p_{ij})$ such that whenever there is an edge $a \to b$ in $\Gamma$, probability $p_{ab}$ is strictly positive.
	
	We say that the sequence of random variables $X_0,X_1,\ldots$ \emph{converges in distribution} to a random variable $X$ (write $X_n \xrightarrow{d} X$) if $\forall g \in G \: \lim_{n \to \infty} Pr(X_n = g) = Pr(X = g)$. The following theorem, which is a direct consequence of the general theory of Markov chains, provides a complete description of the asymptotic behavior of the sequence $\{X_n\}_{n = 0}^{\infty}$ under the assumption that $P$ is doubly stochastic.
	
	\begin{thm}\label{thm-limit}
		For the random walk $\textbf{RW}(G,\: P)$ with doubly stochastic matrix $P$ the following statements hold:
		\begin{gather*}
			\small \textit{If there is a relator of odd length in } G, \textit{ then } X_n \xrightarrow{d} X, \; where \; X \; is \; uniform \; in \; G\\
			\small\textit{Otherwise},X_{2n} \xrightarrow{d} X_{even}, where \; X_{even} \; is \; uniform \; in \; \{g \in G \;|\; \ell(g) \; is \; even\},\\
			\small X_{2n+1} \xrightarrow{d} X_{odd} \; where \; X_{odd} \; is \; uniform \; in \; \{g \in G \;|\; \ell(g) \; is \; odd\}.
		\end{gather*}
	\end{thm}
	
	\begin{cor}\label{cor-main}
		For the random walk \textbf{RW}$(G,\: P)$ with doubly stochastic matrix $P$ and any function $f$ with the domain $G$, the following statements hold:
		\[
		\textit{If there is a relator of odd length in } G, \textit{ then} \lim_{n \to \infty} \EX[f(X_n)] = \frac{1}{|G|} \sum_{g \in G} f(g).
		\]
		\[
		\textit{Otherwise, } \lim_{n \to \infty} \EX[f(X_{2n})] = \frac{2}{|G|} \sum_{g \in G, 2 \mid \ell(g)} f(g) \textit{ and } \lim_{n \to \infty} \EX[f(X_{2n+1})] = \frac{2}{|G|} \sum_{g \in G, 2 \nmid \ell(g)} f(g).
		\]
	\end{cor}
	\begin{proof}
		By Theorem \ref{thm-limit} we know that $X_n \xrightarrow{d} X$, so it follows that $f(X_n) \xrightarrow{d} f(X)$ and $\lim_{n \to \infty} \EX[f(X_n)] = \EX[f(X)] = \frac{1}{|G|} \sum_{g \in G} f(g)$. The other case is done similarly.
	\end{proof}
	
	Further we will often use this corollary for $f= \ell_G$. Now we  compute the limiting length expectation for several examples of groups.
	
	\begin{prp}\label{prp-cyclic}
		For the random walk \textbf{RW}($\mathbb{Z}_{m}, P$) with doubly stochastic matrix $P$ the following statements hold:
		\[
		m \equiv 1,3 \textnormal{ (mod 4)}: \lim_{n \to \infty} \EX \lbrack \ell(X_n) \rbrack = \frac{m}{4} - \frac{1}{4m}.
		\]
		\[
		m \equiv 2 \textnormal{ (mod 4)}: \lim_{n \to \infty} \EX \lbrack \ell(X_{2n}) \rbrack = \frac{m}{4} - \frac{1}{m} \; \textnormal{ and } \lim_{n \to \infty} \EX \lbrack \ell(X_{2n+1}) \rbrack = \frac{m}{4} + \frac{1}{m}.
		\]
		\[
		m \equiv 0 \textnormal{ (mod 4)}: \lim_{n \to \infty} \EX \lbrack \ell(X_n) \rbrack = \frac{m}{4}.
		\]
	\end{prp}
	\begin{proof}
		We will repeatedly use Corollary \ref{cor-main} throughout the proof. Let us consider three cases:
		\begin{itemize}
			\item $m \equiv 1,3\;(\textnormal{mod}\;4)$:
			\[
			\lim_{n \to \infty} \EX[\ell(X_n)] = \frac{\sum_{g \in G} \ell(g)}{m} = \frac{2\left(1 + \ldots + \frac{m-1}{2}\right)}{m} = \frac{m}{4} - \frac{1}{4m}.
			\]            
			\item $m \equiv 2\;(\textnormal{mod}\;4)$:
			\[
			\lim_{n \to \infty} \EX[\ell(X_{2n+1})] = \frac{2}{m} \sum_{g \in G,\: 2 \nmid \ell(g)} \ell(g) = \frac{m}{4} + \frac{1}{m}.
			\]
			\[
			\lim_{n \to \infty} \EX[\ell(X_{2n})] = \frac{2}{m} \sum_{g \in G,\: 2 \mid \ell(g)} \ell(g) = \frac{2}{m} \cdot \left( 2 \left(2 + 4 + \ldots + \frac{m}{2}-1\right) \right) = \frac{m}{4} - \frac{1}{m}.
			\]
			\item $m \equiv 0\;(\textnormal{mod}\;4)$:
			\[
			\lim_{n \to \infty} \EX[\ell(X_{2n+1})] = \frac{2}{m} \sum_{g \in G, 2 \nmid \ell(g)} \ell(g) = \frac{2}{m} \cdot 2 \left(1 + 3 + \ldots + \frac{m}{2} - 1\right) = \frac{m}{4}.
			\]
			\[
			\lim_{n \to \infty} \EX[\ell(X_{2n})] = \frac{2}{m} \sum_{g \in G, 2 \mid \ell(g)} \ell(g) = \frac{2}{m} \left( 2 \left(2 + 4 + \ldots + \frac{m}{2}\right) - \frac{m}{2} \right) = \frac{m}{4}.
			\]
			Hence these two limits coincide and we arrive at 
			\[
			\lim_{n \to \infty} \EX[\ell(X_{2n})] = \lim_{n \to \infty} \EX[\ell(X_{2n+1})] = \frac{m}{4}.
			\]
		\end{itemize}
	\end{proof}
	Given groups $G$ and $H$, denote by $G \times H$ their direct product. It is not hard to see that if $G \cong \langle S_G \; | \; R_G \rangle$ and $H \cong \langle S_H \; | \; R_H \rangle$, then $G \times H \cong \langle S_G \times S_H \; | \; R_G \cup R_H \cup R_{C} \rangle$, where $R_C$ is a set of relations specifying that elements of $S_G$ and $S_H$ commute. So when we further consider a random walk on a direct product of groups, this presentation of a group is meant. Now we will give the formula for the limiting length expectation for the direct product of groups in the following proposition.
	\begin{prp}\label{prp-prod}
		For two non-trivial groups $G_1,G_2$ and the random walk $\textbf{RW}(G_1 \times G_2, P)$ with doubly stochastic matrix $P$ the following holds: \begin{gather*}
			\lim_{n\to \infty} \EX[\ell_{G_1 \times G_2}(X_n)] = \frac{\sum_{g_1 \in G_1} \ell_{G_1}(g_1)}{|G_1|} + \frac{\sum_{g_2 \in G_2} \ell_{G_2}(g_2)}{|G_2|}.
		\end{gather*}
	\end{prp}
	\begin{proof}
		We begin by observing that $\ell_{G_1 \times G_2}((g_1,g_2)) = \ell_{G_1}(g_1) + \ell_{G_2}(g_2)$. If at least one of $G_1$ and $G_2$ has a relator of odd length, then by Corollary \ref{cor-main} we have
		\begin{gather*}
			\lim_{n\to \infty} \EX[\ell_{G_1 \times G_2}(X_n)] =\frac{\sum_{g \in G_1 \times G_2} \ell_{G_1 \times G_2}(g)}{|G_1 \times G_2|} = \frac{\sum_{g_1 \in G_1, g_2 \in G_2} \ell_{G_1}(g_1) + \ell_{G_2}(g_2)}{|G_1| |G_2|} = \\ \frac{|G_2| \sum_{g_1 \in G_1} \ell_{G_1}(g_1) + |G_1| \sum_{g_2 \in G_2} \ell_{G_2}(g_2)}{|G_1||G_2|} = \frac{\sum_{g_1 \in G_1} \ell_{G_1}(g_1)}{|G_1|} + \frac{\sum_{g_2 \in G_2} \ell_{G_2}(g_2)}{|G_2|}.
		\end{gather*}
		
		Otherwise, all relators of $G_1$ and $G_2$ have even length. Notice that in this case the Cayley graph of $G_1$ is regular and does not have odd cycles, therefore it is also bipartite. However, it is not hard to see that this implies that there is the same number of even and odd elements in $G_1$. The same argument applies to $G_2$. Hence, the sum of lengths of even elements is equal to 
		\begin{gather*}
			\sum_{2 \mid \ell_{G_1}(g_1)} \left(\frac{|G_2|}{2}\ell_{G_1}(g_1) + \sum_{2 \mid \ell_{G_2}(g_2)} \ell_{G_2}(g_2) \right) = \frac{|G_2|}{2}\sum_{2 \mid \ell_{G_1}(g_1)} \ell_{G_1}(g_1) + \frac{|G_1|}{2}\sum_{2 \mid \ell_{G_2}(g_2)} \ell_{G_2}(g_2). 
		\end{gather*}
		An analogous expression can be obtained for the sum of lengths of odd elements.
		\begin{gather*}
			\frac{|G_2|}{2}\sum_{2 \nmid \ell_{G_1}(g_1)} \ell_{G_1}(g_1) + \frac{|G_1|}{2}\sum_{2 \nmid \ell_{G_2}(g_2)} \ell_{G_2}(g_2).
		\end{gather*}
		Therefore, if we sum these two expressions and take into account the factor of $\frac{2}{|G_1 \times G_2|} = \frac{2}{|G_1||G_2|}$, the proposition follows from Corollary \ref{cor-main}.
	\end{proof}
	
	\begin{cor}\label{cor-abelian}
		For the random walk \textbf{RW}($\mathbb{Z}_{n_1} \times \cdots \times \mathbb{Z}_{n_k}, P)$ with $k,n_1,\ldots,n_k > 1$ and doubly stochastic matrix $P$ we have the following:
		\[
		\lim_{n \to \infty} \EX[\ell(X_n)] = \sum_{i=1}^{k} \frac{n_i}{4} - \sum_{2 \nmid n_i} \frac{1}{4n_i}.
		\]
	\end{cor}
	\begin{proof}
		It is not difficult to find the average length in the Cayley graph of $\mathbb{Z}_n$. If $n$ is odd, then from Proposition \ref{prp-cyclic} it is equal to $\frac{n}{4} - \frac{1}{4n}$. If $n$ is even, then it is $\frac{2 \cdot(0 + 1 + \ldots + \frac{n}{2}) - \frac{n}{2}}{n} = \frac{n}{4}$. Thus, from Proposition \ref{prp-prod} the result follows. 
	\end{proof}
	
	Applying the fact that each finite abelian group is isomorphic to a direct product of cyclic groups, we can find the limiting length expectation in a certain presentation for any finite abelian group using Corollary \ref{cor-abelian}.

	\section{Coxeter groups, reflection groups and invariant polynomials}
	
	We begin by outlining some terminology and facts regarding finite Coxeter groups and reflection groups, which will be needed later to construct a random walk on a braid group. More details can be found in \cite{humphreys1992reflection}. 
	
	Coxeter group $W$ is a group with presentation $\langle \: s_1, \ldots, s_n \: | \: (s_is_j)^{m(i,j)} = 1 \: \rangle$ where $m(i, i) = 1 \: \forall i = \overline{1,n}$ and $m(i, j) = m(j, i) \ge 2$ is an integer or $\infty$ for $i \neq j$. 
	
	The Artin-Tits braid group of a Coxeter group $W$ is a group with generators $\sigma_1, \ldots, \sigma_n$ subject to the relations:
	$$\underbrace{\sigma_i \sigma_j \ldots}_{m(i,j)} = \underbrace{\sigma_j \sigma_i \ldots}_{m(i,j)}.$$
	
	Let $V$ be a Euclidean space over a field $\mathbb{K}$ of characteristics 0. Define a reflection $s_\alpha$ as a nonidentical operator of the Euclidean space, which fixes pointwise the hyperplane $H_\alpha$ orthogonal to the vector $\alpha$. Groups generated by a finite set of reflections are called finite reflection groups. It occurs that finite Coxeter groups are precisely finite reflection groups.
	
	Let $G$ be a finite subgroup of $GL(V)$. Denote by $S$ the symmetric algebra on dual space $V^{*}$. Fixing a basis in $V$, $S$ can be performed as an algebra $\mathbb{K}[x_1,\ldots,x_n]$ where $x_1,\ldots,x_n$ are coordinate functions. There is a natural action of $G$ on $S$ $(g \cdot f)(v) = f(g^{-1}v)$, where $g \in G, \: v \in V, \: f \in S$. In addition, $f$ is said to be $G$-invariant if $g \cdot f = f \; \forall g \in G$. 
	\begin{thm}(Chevalley)\label{thm-chev}
		Every subalgebra of $\mathbb{K}\left[x_1,\ldots,x_n\right]$ consisting of $W$-invariant polynomials is generated as an $\mathbb{K}$-algebra by $n$ algebraically independent homogeneous elements of positive degree (together with 1). 
	\end{thm}
	Polynomials from the above theorem are called \emph{basic invariant polynomials}, and denote their degrees by $d_1, \ldots, d_n$. Although there can be many sets of generators, their degrees turn out to be unique up to reordering.
	
	We also provide a table of degrees of basic invariant polynomials for irreducible Coxeter groups. In addition, any finite Coxeter group $W$ can be represented as a direct product of such groups.
	\begin{center}
		
		\begin{tabular}{ | C{2.5cm} | C{4cm} | } 
			\hline
			$W$ & $d_1,\ldots,d_n$ \\
			\hline
			$A_n$ & $2,3,\ldots,n+1$ \\ 
			\hline
			$B_n/C_n (n \ge 2)$ & $2,4,6,\ldots,2n$ \\ 
			\hline
			$D_n (n \ge 4)$ & $n; 2,4,6,\ldots,2n-2$ \\ 
			\hline
			$E_6$ & $2, 5, 6, 8, 9, 12$ \\ 
			\hline
			$E_7$ & $2, 6, 8, 10,
			12, 14, 18$ \\ 
			\hline
			$E_8$ & $2, 8, 12, 14,
			18, 20, 24, 30$ \\ 
			\hline
			$F_4$ & $2, 6, 8, 12$ \\ 
			\hline
			$G_2$ & $2, 6$ \\ 
			\hline
			$H_3$ & $2, 6, 10$ \\ 
			\hline
			$H_4$ & $2, 12, 20, 30$ \\ 
			\hline
			$I_2(m)(m \ge 3)$ & $2, m$ \\ 
			\hline
		\end{tabular}
	\end{center}
	
	Define a map $w \mapsto \widetilde{w}$ from $W$ to the Artin-Tits braid group of $W$ as follows: if $w = s_{i_1} \ldots s_{i_r}$ is a reduced expression of $w$ in $W$, then $\widetilde{w} := \sigma_{i_1} \ldots \sigma_{i_r}$. This mapping is well-defined by Matsumoto's Theorem stated below. 
	\begin{thm}(Matsumoto)\label{thm-mats}
		Any two reduced expressions of $w \in W$ are connected by a sequence of braid moves.
	\end{thm}
	Although $w \mapsto \widetilde{w}$ is not a group homomorphism, it still preserves the length function $\ell_{{W}}(w) = \ell_{\widetilde{W}}(\widetilde{w})$. Now using this map, we will construct a random walk on a braid group by first considering a random walk on a Coxeter group $W$ and then lifting it into $\widetilde{W}$. 
	
	In the next proposition, we will compute the limiting expectation of the braid length in our model.
	\begin{prp}\label{prp-coxt}
		If $W \neq W_{A_1}, W_{{I}_2(m)}$ for odd $m$, then for the random walk \textbf{RW}($W, P$) with double stochastic matrix $P$ we have the following:
		\[
		\lim_{n \to \infty} \EX[\ell_{\widetilde{W}}({\widetilde{X}}_n)] = \frac{Ref\:W}{2} 
		\]
		where $Ref\;W$ is the reflection number of $W$, which is equal to the length of the longest element in $W$.
	\end{prp}
	\begin{proof}
		First, recall that $\ell_{\widetilde{W}}(\tilde{\beta})=\ell_W(\beta)$. It is well known that the generating function of $\ell(w)$ is given by its Poincaré polynomial which admits the following factorization:
		\[
		P_W(t) := \sum_{w \in W} t^{\ell(w)} = \prod_{i=1}^{n} \frac{t^{d_i}-1}{t-1}.
		\] Differentiate it with respect to $t$ in order to obtain
		\[
		\sum_{w \in W} \ell(w) \cdot t^{\ell(w)-1} = \sum_{i=1}^{n} (1+\ldots+t^{d_1-1}) \ldots (1+\ldots+t^{d_n-1}) \cdot \frac{(1 + 2t + \ldots + (d_i-1)t^{d_i-2})}{(1+\ldots+t^{d_i-1})}.
		\]
		Evaluating it at $t=1$ using the identities $d_1 \cdot \ldots \cdot d_n = P_W(1) = |W|$, $d_1 + \ldots + d_n - n = deg \; P_W = Ref\:W$, gives us the sum of the lengths of all elements in $W$
		\[
		\sum_{w \in W}\ell(w) = P_W^{\prime}(1) = \frac{|W|}{2} \cdot \sum_{i=1}^{n} (d_i -1) = \frac{Ref \: W}{2} \cdot|W|.
		\]    
		From the classification of finite Coxeter groups, we know that $W$ is a direct product of irreducible Coxeter groups that are listed in the table above. Additionally, a multiset of degrees of basic $W$-invariant polynomials is a union of degree multisets of its individual irreducible components (in the sense that multiplicities of respective elements are added). However, $W \neq W_{A_1},W_{{I}_2(m)}$, so it can be seen from the table above that $W$ has at least two basic invariant polynomials of even degree. Note that the case of $W_{A_2}$ is already covered since Coxeter groups $W_{A_2}$ and $W_{I_2(3)}$ have the same presentations. Denote their degrees $d_i$ and $d_j$, then each summand in $P_W^{\prime}(-1)$ contains either $1+\ldots+q^{d_i - 1}$ or $1+\ldots+q^{d_j - 1}$ both of which vanish at $q=-1$, so we get $\sum_{w \in W} (-1)^{\ell(w)} \ell(w) = 0$ from which we deduce that 
		\[
		\sum_{w \in W,\:2 \mid \ell(w)} \ell(w) = \sum_{w \in W,\:2 \nmid \ell(w)} \ell(w) = |W| \cdot \frac{Ref \: W}{4}.
		\]
		Since all relations in the presentation of $W$ have even length by Corollary \ref{cor-main} we obtain the desired result.
	\end{proof}
	
	Now notice that in the exceptional cases we do actually get two distinct limits. The case of $W_{A_1}$ is trivial and observe that the case of $ W_{{I}_2(m)}$ with odd $m$ was already considered in Proposition \ref{prp-cyclic} since the Cayley graph of $W_{I_2(m)}$ is a $2m$-cycle, which also gives us two distinct limits.
	\[
	\lim_{n \to \infty} \EX \lbrack \ell(X_{2n}) \rbrack = \frac{m}{2} - \frac{1}{2m} \; \textnormal{ and } \lim_{n \to \infty} \EX \lbrack \ell(X_{2n+1}) \rbrack = \frac{m}{2} + \frac{1}{2m}.
	\]
	
	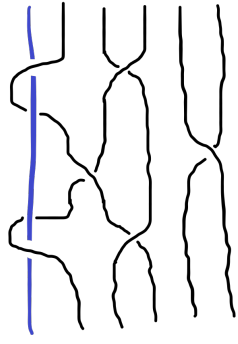
\begin{figure}[h!]
		\centering
		\begin{tikzpicture}
            \draw[color=blue, line width=1pt] (0,0)--(0,-1.9);
            \draw[black, line width=1pt] plot [smooth] coordinates {(1,0) (1,-1.5) (-.3, -2.25) (1,-3) (1,-4.5)};
            \filldraw[white] (0,-2.5) circle (3pt);
            \draw[color=blue, line width=1pt] (0,-2.15)--(0, -4.5);
        
            \braid[line width=1pt, nudge factor=-0.05, control factor=0.75, number of strands=4] (braid) at (2,0) s_1 s_3 s_2^{-1} s_3^{-1};
        \end{tikzpicture}
		\caption{Example of a torus braid}
		\label{fig:enter-label}
	\end{figure}
	
	\emph{Remark}. There is a geometrical realization of the random walk in type $B$ obtained by lifting the braids into the so-called "torus braid group" $Br_{1,n}$. Torus braids are constructed in a similar way to usual braids. We just consider a braid on $n+1$ strands where one strand is fixed as shown in Figure 4. Then define $Br_{1,n}$ to be the group of equivalence classes of torus braids under ambient isotopy with concatenation of braids as the group operation. In \cite{lambropoulou1993study} it is proved that $Br_{1,n}$ is isomorphic to $\widetilde{W}_{{B_n}}$.
	
	Recall that the symmetric group $S_n$ is a Coxeter group of type $A_{n-1}$. We will now find the limiting expectation for the number of connected components of random links obtained by closing braids in the random walk on $Br_n$.
	
	\begin{prp}
		For the random walk $\textbf{RW}(S_n, P)$ with doubly stochastic matrix $P$ the following holds:
		\begin{gather*}
			\lim_{N \to \infty} \EX[c(\overline{\widetilde{X}_{2N}})] = H_n - \frac{(-1)^{n}}{n(n-1)}, \\
			\lim_{N \to \infty} \EX[c(\overline{\widetilde{X}_{2N+1}})] = H_n + \frac{(-1)^{n}}{n(n-1)}
		\end{gather*}
		where $c(L)$ denotes the number of connected components of a link $L$ and $H_n = \sum_{i=1}^n \frac{1}{i}$.
	\end{prp}
	\begin{proof}
		Let $c(w)$ denote the number of cycles in the cycle decomposition of $w \in S_n$, then $c(w) = c(\overline{\widetilde{w}})$ since each cycle in the cycle decomposition of $w \in S_{n}$ after closure will become a connected component of an obtained link. It is well known that the generating function for $c(w)$ is 
		\begin{gather*}
			F(x) =\sum_{w \in S_n} x^{c(w)} = x(x+1) \ldots (x+n-1).
		\end{gather*}
		If we differentiate it and evaluate at $x=1$, we will obtain
		\begin{gather*}
			F'(x) = x(x+1)\ldots (x+n-1) \cdot \left(\frac{1}{x} + \ldots + \frac{1}{x+n-1}\right), \\
			F'(1) = \sum_{w \in S_n} c(w) = n! H_n.
		\end{gather*}
		Now evaluation at $x=-1$ and the fact that $\ell(w)$ and $n-c(w)$ have the same parity, give us that
		\begin{gather*}
			x^{n+1}F'(x) = \sum_{w \in S_n} c(w)x^{n+c(w)}, \\
			\sum_{w \in S_n} c(w) (-1)^{\ell(w)} =  \sum_{w \in S_n} c(w) (-1)^{n-c(w)} = \sum_{w \in S_n} c(w) (-1)^{n+c(w)} = (-1)^{n+1}F'(-1) = (-1)^{n}(n-2)!.
		\end{gather*}
		Now putting these two results together and using Corollary \ref{cor-main}, we will get
		\begin{gather*}
			\lim_{N \to \infty} \EX[c(\overline{\widetilde{X}_{2N}})] = \frac{2}{n!} \frac{n!H_n + (-1)^n (n-2)!}{2} = H_n + \frac{(-1)^{n}}{n(n-1)}, \\
			\lim_{N \to \infty} \EX[c(\overline{\widetilde{X}_{2N+1}})] =\frac{2}{n!} \frac{n!H_n - (-1)^n (n-2)!}{2} = H_n - \frac{(-1)^{n}}{n(n-1)}
		\end{gather*} which finishes the proof.
		
	\end{proof}
	
	Observe that $H_n \sim \log n$, and thus both expected numbers of connected components are asymptotically equivalent to $\log n$ as $n\to \infty$.
	
	\section{Large deviations}
	
	Let $\{X_i\}_{i=1}^{\infty}$ be a sequence of uniform i.i.d. random variables in $\{\sigma \in S_n \;|\; \ell(\sigma)\;is\;even\}$. Given $\widetilde{X}_1 = \beta_1, \ldots, \widetilde{X}_N = \beta_N$ we will now shift indices for each separate braid and compose them into a "block-diagonal form", so that any two of them commute.   
	
	\begin{figure}[!h]
		\centering
		\begin{tikzpicture}[scale=.8]
			\braid[line width=1pt, nudge factor=-0.05, control factor=0.75] (braid) at (0,0) s_1 s_2^{-1} s_{4} s_4 s_{6};
			\draw[draw=gray] (-0.3,0.3) rectangle (2.3,-2.5);
			\draw[draw=gray] (2.7,-2.5) rectangle (4.3,-4.2);
			\draw[draw=gray] (4.7,-4.3) rectangle (6.3,-5.8);
			\node at (-.8,0) {$\widetilde{X}_1$};
			\node at (2.35, -3.35) {$\widetilde{X}_2$};
			\node at (6.7, -5.05) {$\widetilde{X}_3$};
		\end{tikzpicture}
		\caption{Example for $\widetilde{X}_1 = \sigma_{1} \sigma_2^{-1}, \widetilde{X}_2 = \sigma_1^2, \widetilde{X}_3 = \sigma_1$.}
		\label{fig:enter-label}
	\end{figure}
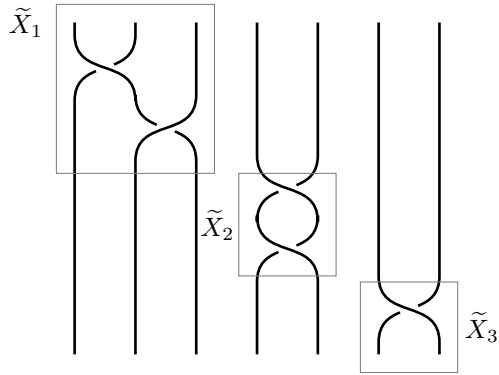
	
	The random variable describing the length function of the obtained braid composition is 
	\[
	L_N := \sum_{i=1}^{N} \ell_{Br}(\widetilde{X}_i).
	\]
	Next, we will show that the probability of the event $L_N = 2 Nx $ where $Nx$ is an integer and $x \in \left[0;\: M = \lfloor \frac{n(n-1)}{4} \rfloor \right]$  is exponentially small. 
	\begin{prp}
		\[
		\lim_{N \to \infty} \frac{\log {Pr}(L_N = 2 Nx )}{N} = -I(x) \textnormal{ where} 
		\]
		\[
		I(x) = \begin{cases}
			\left(1 - \frac{1}{M+1 - x}\right)^{M+1} + (M-x)\log(M-x) - (M+1-x)\log(M+1-x) + \log\left(\frac{n!}{2}\right), \textnormal{ for } x \in [0;\:M] \\
			
			+\infty, \textnormal{ otherwise}.
		\end{cases}
		\]
	\end{prp}
	\begin{proof}
		Let $\kappa(n,j,k)$ denote the number of ordered integer tuples $(x_1,\ldots,x_k)$ such that
		\[
		x_1 + \ldots + x_k = n \text{ and } 0 \le x_i < j.
		\]
		From \cite{419766} we know the asymptotics of the number of restricted compositions $\kappa( Nx , M+1, N)$
		\[
		\kappa( Nx, M+1, N) \sim \frac{1}{\sqrt{2 \pi N}} \cdot \frac{(M+1-x)^{N(M+1-x)+\frac{1}{2}}}{(M-x)^{N(M-x)+\frac{3}{2}}} \cdot \exp\left({-N \: \left(\frac{M-x}{M+1-x}\right)^{M+1}}\right).
		\]
		The number of solutions to the equation 
		\[
		\frac{y_1}{2} + \ldots + \frac{y_N}{2} =  Nx  \text{ with } 0 \le \frac{y_i}{2} \le M \text{ and } y_i \text{ even for } i = \overline{1,N}.
		\]
		is exactly given by $\kappa( Nx , M+1, N)$.
		Since $X_1, \ldots, X_N$ are uniform i.i.d. random variables  in $\{\sigma \in S_n \;|\; \ell(\sigma)\;is\;even\}$, then $Pr(\widetilde{X}_1 = x_1, \ldots, \widetilde{X}_{N}=x_N) = \left(\frac{2}{n!} \right)^N$. Thus, it is not difficult to verify the asymptotics
		\[
		{Pr}(L_N = 2 Nx ) = \kappa( Nx , M+1, N) \left(\frac{2}{n!}\right)^{N} = \exp\left(-NI(x) + O(\log N)\right)
		\] from which the desired result immediately follows.
	\end{proof}
	
	\printbibliography
	
\end{document}